\documentclass{article}

\usepackage{arxiv}

\usepackage[english]{babel}
\usepackage{biblatex}
\usepackage{amsmath}
\usepackage[pdftex]{hyperref}
\usepackage[utf8]{inputenc} % allow utf-8 input
\usepackage[T1]{fontenc}    % use 8-bit T1 fonts
\usepackage{hyperref}       % hyperlinks
\usepackage{url}            % simple URL typesetting
\usepackage{booktabs}       % professional-quality tables
\usepackage{amsfonts}       % blackboard math symbols
\usepackage{nicefrac}       % compact symbols for 1/2, etc.
\usepackage{microtype}      % microtypography
\usepackage{lipsum}
\usepackage{graphicx}
\graphicspath{ {./images/} }

\title{A formula to solve sextic degree equation}

\author{
 Rodrigo José Martinelli Biglia Andrade \\
  Administration Department\\
  Methodist University of Piracicaba\\
  Piracicaba, SP, Brazil \\
  \texttt{contato@pontoabc.com} \\
}

\begin{document}
\maketitle
\begin{abstract}
According to the Abel-Ruffini theorem, equations of degree
equal to or greater than 5 cannot, in most cases, be solved by
radicals. Due of this theorem we will present a formula that solves specific cases of sixth degree equations using Martinelli’s
polynomial as a base. To better understand how this formula
works, we will solve a sixth degree equation as an example.
We will also see that all sixth degree equations that meet the
coefficient criterion have a resolvent of fifth degree that can be splitted into a second degree and a third degree
equation. Throughout the paper we will see a demonstration of
the ratio of the coefficients of a sixth degree equation that can
be solved with the formula that will be presented this paper
\end{abstract}

% keywords can be removed

\section{Introduction}
First, this paper aims to contribute to the academic content on
mathematics and focuses on professionals in the mathematic area who wish to increase the range of content they may be teaching in class to their students. Second, the purpose of this paper is also to present in a didactic way of a formula (the Milanez's formula) which it is able to solve specific cases of the sixth degree equation. This formula is deducted by the Martinelli's polynomial and consequently we will have a fifth degree resolvent of the sixth degree equation. Every fifth degree equation that is splittable into a degree 2 and a degree 3 is the resolvent of a sixth degree equation. The demonstration of the relation of the sixth degree equation will be showed a posteriori, just in cases where the coefficients are consistent with the relation that will be presented and proved mathematically it is possible to solve a sixth degree equation where it has a splittable fifth degree resolvent.
Using an example, we will have a better understanding of solving a sixth degree equation according to Milanez's relation.
All equations that obey the Milanez's relation are solvable by radicals where the roots of the sixth degree polynomial are the sum of the roots of a polynomial of degree 2 and with a polynomial of degree 3. However, the equation of sixth degree has a resolvent of fifth degree where it is not necessary to use any transformation to eliminates the fourth degree term of the fifth degree equation because all resolvent will have the fourth degree term suppressed.
It is convenient to demonstrate and explain every detail of the Martinelli's polynomial so that we can give greater understanding of the development of Milanez's formula.
But, it is also interesting for the reader to know that any of the roots of Martinelli's polynomial has the unique utility of be used to split an equation of the fifth degree into two of degree minor as a degree 2 and a degree 3. So Martinelli's equality is also useful for solving fifth degree equations. In this paper we will present a very simple way to use this feature to be able to deduce a formula for specific cases of equation of the sixth degree where its roots can be represented by means of radicals.
This paper also commits itself to a friendly and didactic content for the reader, making the inviting and interesting reading for any student and teacher in the mathematic field.

\section{About Galois theory and Abel-Ruffini theorem}
\label{sec:headings}
We will see some concepts about Galois' theory and Abel-Ruffini's theorem.
The theorem, in a succinct way, consists of the proof that there is no "closed" formula for all equations of degree greater or equal to 5. That is, there is no way to algebraically arrive at a formula for those equations greater or equal to 5. As Zoladek (2000) reported, "A general algebraic equation of degree $\geq$ 5 cannot be solved in radicals. This means that there does not exist any formula
which would express the roots of such equation as functions of the coefficients by means of the algebraic operations and roots of natural degrees."  (p. 254).

Thus, Galois theory comes to confirm, through the study of the relationship of the roots that it is impossible to solve, algebraically, equations of degree equal to or greater than 5. In its most basic form, the theorem states that given a
E / F field extension which is finite and Galois, there is a one-to-one correspondence between its intermediate fields and subgroups of its Galois group. (Intermediate fields are K fields that satisfy $F \subseteq K \subseteq E$; they are also called sub-extensions of E / F.).

So, according to Abel-Ruffini's theorem and Galois theory, there are cases where equations of degree 5 or higher are solvable by radicals or algebraically.

\section{Proof of Martinelli's polynomial}

To demonstrate the Martinelli's Polynomial it is necessary to consider the roots of any fifth degree equation.
Each root must be associated with a single other root. We have then that the combination of the roots of a fifth degree equation is ten (due that the polynomial is tenth degree). Let us then consider the following equations:

\begin{equation}
(x^2-kx+n)(x^3+kx^2+kx+m)=0.
\end{equation}
\begin{equation}
(x^2-kx+n)(x^3+kx^2+lx+m)=0.
\end{equation}

Equation (1) have a relation with equation (2). If we expand both equations (1) and (2), each corresponding term can be splittable. Let put those terms into a table:

\def\arraystretch{2.0}%
\begin{table}[htb]
\centering
\caption{Equation (1) expanded}
\begin{tabular}{|c|c|c|c|}
\hline
$x^3$ & $x^2$ & $x$ & Term Independent \\
\hline
$k+n-k^2=C_2$ & $-k^2+kn+m=D_2$ & $kn-km=E_2$ & $ mn=F$ \\
\hline
\end{tabular}
\end{table}
\def\arraystretch{2.0}%
\begin{table}[htb]
\centering
\caption{Equation (2) expanded}
\begin{tabular}{|c|c|c|c|}
\hline
$x^3$ & $x^2$ & $x$ & Term Independent \\
\hline
$l+n-k^2=C$ & $nk-lk+m=D$ & $nl-km=E$ & $ mn=F$ \\
\hline
\end{tabular}
\end{table}

Andrade (2019) consider the follow equation to better understand the proof of Martinelli's Polynomial:

\begin{equation}
x^2+3x+2=0
\end{equation}

If we substitute x for the sum of the roots of equation (3) we have 2=0, that would be absurd, but with that idea, we can create a second degree equation using the terms $C, C_2 , D, D_2, E$ e $E_2$ as shown in Tables 1 and 2 above. Matching the equation that will be created with $(-k^2+n+k-C)n$ we will have, on the right side of the equation, the same thing $E_2- E$, because $E_2$ is the same as $kn-km$ and $-E$ is the same as $-nl+km$. The sum of $E_2-E$ will result $kn-nl$ then $l=C-n+k^2$. So we have $kn -n(C-n+k^2)$ that is the same as $(-k^2+n+k-C)n$.

To create a second degree equation where one of the solutions of the equation will be the sum of two solutions of a fifth degree equation we must follow this logic: We have the general form of the second degree equation that is
$ax^2 + bx + c = 0$, the coefficients of the second degree equation that will be created will be $a = C_2- C$, $b = D_2- D$, $c=E_2- E$ and on the right side of the equation we have to add $(-k^2+n+k-C)n$ because on the left side of the equation will remain $E_2- E$.

This idea will give us the possibility to know n from equation (2). The goal is to form a tenth degree equation with the unknown k. So $a = C_2 - C$ e $b=D_2-D$ according to tables 1 and 2, it follows that $a=k+n-k^2-C$, $b=-k^2+kn+m-D$, $c=nk -km-E$ and the right side of the equation will be $(-k^2+n+k-C)n$. Thus, making the substitutions in $ak^2+bk+c=0$ we will have:

\begin{align*}
(-k^2+n+k-C)k^2 + (-k^2+nk+m-D)k + nk-km-E &= (-k^2+n+k-C)n.\\
(-k^4+nk^2+k^3-Ck^2)+(-k^3+nk^2+km- Dk)+ nk-km- E &= (-k^2+n+k-C)n. \\
-k^4+2nk^2-Ck^2-Dk+ nk- E &= -nk^2+n^2+nk-nC.\\
-k^4-Ck^2-Dk-E&=-3nk^2+n^2-nC.\\
\end{align*}

Putting n on the left side of the equation, we have:

\begin{equation}
n = \frac{k^4+Ck^2+Dk+E}{3k^2-n+C}.
\end{equation}

Since k is the sum of two roots of a fifth degree equation, C, D, E and F the coefficients of a fifth degree equation and n is the product of the roots of a second degree equation, so to arrive at Martinelli's polynomial, we have to replace the variable n with the use of an algebraic manipulation. This algebraic manipulation consists of taking the equality referring to the term D in table (2). Thus:

\begin{equation}
nk-(k^2-n+C)k+\frac{F}{n}=D.
\end{equation}

\begin{equation}
n^2=\frac{nk^3+Cnk-F+Dn}{2k}.
\end{equation}

Now just replace $ n^2$ in equation (4) and get n:

\begin{equation*}
n = \frac{k^4+Ck^2+Dk+E}{3k^2-n+C}.
\end{equation*}
\begin{equation*}
-n^2+3nk^2+Cn=k^4+Ck^2+Dk+E.
\end{equation*}
\begin{equation*}
\frac{-nk^3-Cnk+F-Dn}{2k} + 3nk^2 + Cn = k^4+Ck^2+Dk+E.
\end{equation*}
\begin{equation}
n = \frac{2(k^5+Ck^3+Dk^2+Ek)-F}{5k^3+Ck-D}.
\end{equation}

If we have n, then we can create the Martinelli's polynomial that will be very important to Milanez's formula.
In equality (4) we can replace the n on the left and the right side of the equation, like this:

\begin{equation*}
\frac{2(k^5+Ck^3+Dk^2+Ek)-F}{5k^3+Ck-D} = \frac{k^4+Ck^2+Dk+E}{3k^2-\frac{2(k^5+Ck^3+Dk^2+Ek)-F}{5k^3+Ck-D}+C}.
\end{equation*}

Arranging the equation on both sides and equaling 0 we arrive at the Martinelli polynomial:

\begin{eqnarray*}
{(2(k^5+Ck^3+Dk^2+Ek)-F)(13k^5+6Ck^3-5Dk^2+(-2E+C^2)k+F-DC)}
\end{eqnarray*}
\begin{equation}
-(k^4+Ck^2+Dk+E)(5k^3+Ck-D)^2 = 0.
\end{equation}
\section{Milanez's relationship}

When a fifth degree equation is perfectly separable into a degree 2 equation and another equation of degree 3, then we can say that the Martinelli's polynomial is also perfectly separable into two, in this case, one of degree 4 and another of degree 6.

The following list shows what a splittable equation looks like so that it is algebraically represented, that is, by radicals.

\begin{equation}
(x^2-ax+b)(x^3+ax^2+cx+d) = 0. \\
\end{equation} 

Expanding equation (9) we arrive at a fifth degree equation of the form:

\begin{equation}
x^5+(-a^2+c+b)x^3+(ab+d-ac)x^2+(bc-ad)x+bd=0. \\
\end{equation}

Knowing that C = $ (-a^2+c+b)$, D=$(ab+d-ac)$, E=$ (bc-ad)$ and F=$bd$, then we can create a relationship with the Martinelli's polynomial which, if expanded, will have the following terms:

\begin{equation*}
k^{10} + 3Ck^8 + Dk^7 + (3C^2-3E)k^6 + (2DC-11F)k^5 + (C^3-D^2-2CE)k^4 + (DC^2-4DE-4CF)k^3
\end{equation*}
\begin{equation}
 + (7DF-CD^2-4E^2+EC^2)k^2 + (4EF-FC^2-D^3)k-F^2+FDC-D^2E=0.
\end{equation}

If we substitute each letter for the relation of equation (10) then we have a polynomial that can be splittable  into a fourth and a sixth degree. Thus:

\begin{equation*}
(k^4+(a)k^3+(c-a^2)k^2+(-a^3-d)k+(ad-a^2c))
\end{equation*}
\begin{equation*}
(k^6+(-a)k^5+(-a^2+3b+2c)k^4+(a^3-2ab-2ac+2d)k^3+(-a^2b+c^2+3b^2-ad)k^2
\end{equation*}
\begin{equation}
+(-ac^2-ab^2+a^2d+2abc-6bd+2dc)k+(adb-adc-2b^2c+d^2+b^3+bc^2))=0.
\end{equation}

So we can solve specific cases of the sixth degree equation by radicals when the fifth degree resolvent equation is perfectly separable. So, the Milanez's relation are the coefficients of each term of the sixth degree equation that makes up the Martinelli's polynomial when it can be separated into two equations (one of degree 4 and one of degree 6). Thus, we have that the relation to obtain a sixth degree equation solvable by radicals when the resolvent is an equation of fifth degree splittable into one of degree 2 and another of degree 3. So, the Milanez's relation is:

\begin{align*}
k^6\\
-ak^5\\
(-a^2+3b+2c)k^4 \\
(a^3-2ab-2ac+2d)k^3\\
(-a^2b+c^2+3b^2-ad)k^2\\
(-ac^2-ab^2+a^2d+2abc-6bd+2dc)k\\
(adb-adc-2b^2c+d^2+b^3+bc^2)\\
\end{align*}

\section{Milanez's formula}

To deduce Milanez's formula is quite simple. If the Martinelli's polynomial is the sum of two roots of any equation of any fifth degree, then the formula for solving a sixth degree equation where the coefficients are consistent with the given Milanez's relation is:

\begin{equation*}
\begin{split}
x=\sqrt[3]{\bigg({-\frac{a_2^3}{27a_1^3}}+\frac{{a_2}{a_3}}{6a_1^2}-\frac{{a_4}}{2a_1}\bigg)+ \sqrt{\bigg(-\frac{{a_2^3}}{27a_1^3}+\frac{{a_2}{a_3}}{6a_1^2}-\frac{a_4}{2a_1}\bigg)^2+\bigg(\frac{a_3}{3a_1}-\frac{a_2^2}{9a_1^2}\bigg)^3}} + \\ 
\sqrt[3]{\bigg({-\frac{a_2^3}{27a_1^3}}-\frac{{a_2}{a_3}}{6a_1^2}-\frac{{a_4}}{2a_1}\bigg)-\sqrt{\bigg(-\frac{{a_2^3}}{27a_1^3}+\frac{{a_2}{a_3}}{6a_1^2}-\frac{a_4}{2a_1}\bigg)^2+\bigg(\frac{a_3}{3a_1}-\frac{a_2^2}{9a_1^2}\bigg)^3}} + 
\end{split}
\end{equation*}
\begin{equation}
\frac{-{a_5}+\sqrt{a_5^2-4{a_6}}}{2}.
\end{equation}

In equation (13), the coefficients $ a_1, a_2, a_3, a_4, a_5, a_6 $ are from two equations, one from the third and the other from the second degree.

Equation (13) is then the formula that provides a root of the sixth degree equation that can be solved by radicals according to the Milanez's relation. We will see in more details with the resolution of a sixth degree equation.

\subsection{Example: Solution of a sixth degree equation}

As an example, we will solve a sixth degree equation according to Milanez's relation.
The equation is: \begin{equation}x^6+2x^3+21x^2-18x+51=0.\end{equation}

In equation (14) the term a is 0, and based on this we can find the terms b, c and d. For that we will get Milanez's relation and make the comparisons.

If $a=0$ then, according to the relation we have:

\begin{align*}
-a &= 0. \\
3b + 2c &= 0. \\
d &= 1. \\
c^2 + 3b^2 &= 21. \\
-6b + 2c &= - 18. \\
-2b^2c+b^3 + bc^2 &= 50.
\end{align*}
Now we  have a system with two variables that we can solve without problems. It is convenient to choose the simplest equations.
\begin{equation}
3b + 2c = 0.
\end{equation}
\begin{equation}
-6b + 2c = -18.
\end{equation}
If we solve the system of the equation (15) and (16) we will have $b=2$ and $c=-3$. Then we find all the variables in our equation according to the Milanez's relation. So just substitute in equation (9) and we have a fifth degree resolvent equation. Let's see:
\begin{equation*}
(x^2-0x+2)(x^3+0x^2-3x+1) = 0.
\end{equation*}
\begin{equation}
x^5-x^3+x^2-6x+2=0.
\end{equation}

Solving the second degree and third degree equations and adding the roots of the second degree equation with each root of the third degree equation, we will have the roots of the sixth degree equation, that is, we will be using the Milanez's formula to be able to obtain the roots of equation (14). Thus, according to formula (13) one of the roots of equation (14) is:

\begin{equation*}
\begin{split}
x=\sqrt[3]{{-\frac{(0)^3}{27(1)^3}}+\frac{{(0)}{(-3)}}{6(1)^2}-\frac{{1}}{2(1)}\bigg)+ \sqrt{\bigg(-\frac{{(0)^3}}{27(1)^3}+\frac{{(0)}{(-3)}}{6(1)^2}-\frac{(1)}{2(1)}\bigg)^2+\bigg(\frac{(-3)}{3(1)}-\frac{(0)^2}{9(1)^2}\bigg)^3}} + \\ 
\sqrt[3]{\bigg({-\frac{(0)^3}{27(1)^3}}+\frac{{(0)}{(-3)}}{6(1)^2}-\frac{{1}}{2(1)}\bigg)- \sqrt{\bigg(-\frac{{(0)^3}}{27(1)^3}+\frac{{(0)}{(-3)}}{6(1)^2}-\frac{(1)}{2(1)}\bigg)^2+\bigg(\frac{(-3)}{3(1)}-\frac{(0)^2}{9(1)^2}\bigg)^3}} +
\end{split}
\end{equation*}
\begin{equation*}
\frac{-{(0)}+\sqrt{(0)^2-4{(2)}}}{2}.
\end{equation*}

\begin{equation}
x=\sqrt[3]{{-\frac{1}{2}+\sqrt{-\frac{3}{4}}}} + \\
\sqrt[3]{{-\frac{1}{2}-\sqrt{-\frac{3}{4}}}} + \\
\sqrt{2}i.
\end{equation}

Or approximately $1.5320888 .... + 1.4142 ..... i$.
The sum of the roots of the second degree equation with the roots of the third degree equation forms the roots of a sixth degree equation according to the Milanez's relation.
Then the other roots will be:

\begin{equation}
x_2=\sqrt[3]{{-\frac{1}{2}+\sqrt{-\frac{3}{4}}}} +
\sqrt[3]{{-\frac{1}{2}-\sqrt{-\frac{3}{4}}}} -
\sqrt{2}i.
\end{equation}
\begin{equation}
x_3=\frac{1}{2}\bigg(-1+\sqrt{3}i\bigg)\sqrt[3]{{-\frac{1}{2}+\sqrt{-\frac{3}{4}}}} +
\frac{1}{2}\bigg(-1-\sqrt{3}i\bigg)\sqrt[3]{{-\frac{1}{2}-\sqrt{-\frac{3}{4}}}} +
\sqrt{2}i.
\end{equation}
\begin{equation}
x_4=\frac{1}{2}\bigg(-1+\sqrt{3}i\bigg)\sqrt[3]{{-\frac{1}{2}+\sqrt{-\frac{3}{4}}}} +
\frac{1}{2}\bigg(-1-\sqrt{3}i\bigg)\sqrt[3]{{-\frac{1}{2}-\sqrt{-\frac{3}{4}}}} -
\sqrt{2}i.
\end{equation}
\begin{equation}
x_5=\frac{1}{2}\bigg(-1-\sqrt{3}i\bigg)\sqrt[3]{{-\frac{1}{2}+\sqrt{-\frac{3}{4}}}} +
\frac{1}{2}\bigg(-1+\sqrt{3}i\bigg)\sqrt[3]{{-\frac{1}{2}+\sqrt{-\frac{3}{4}}}} +
\sqrt{2}i.
\end{equation}
\begin{equation}
x_6=\frac{1}{2}\bigg(-1-\sqrt{3}i\bigg)\sqrt[3]{{-\frac{1}{2}+\sqrt{-\frac{3}{4}}}} +
\frac{1}{2}\bigg(-1+\sqrt{3}i\bigg)\sqrt[3]{{-\frac{1}{2}-\sqrt{-\frac{3}{4}}}} -
\sqrt{2}i.
\end{equation}

Or if you prefer, we can represent the roots in a trigonometric way, like this:

\begin{align*}
x&=2\cos\bigg(\frac{2\pi}{9}\bigg)+\sqrt{2}i. \\
x_2&=2\cos\bigg(\frac{2\pi}{9}\bigg)-\sqrt{2}i. \\
x_3&=-2\cos\bigg(\frac{\pi}{9}\bigg)+\sqrt{2}i. \\
x_4&=-2\cos\bigg(\frac{\pi}{9}\bigg)-\sqrt{2}i. \\
x_5&=\cos\bigg(\frac{\pi}{9}\bigg)-\sqrt{3}\sin\bigg(\frac{\pi}{9}\bigg) + \sqrt{2}i. \\
x_6&=\cos\bigg(\frac{\pi}{9}\bigg)-\sqrt{3}\sin\bigg(\frac{\pi}{9}\bigg) - \sqrt{2}i.
\end{align*}

In the case of example equation (14) all roots are complex as seen above.

\section{Remark}

As seen, Milanez's relation is fundamental to solving equations of the sixth degree. All sixth degree equations that obey the coefficients according to Milanez's relation can be solved by radicals.
As shown, Milanez's formula is based on the Martinelli's polynomial which has the function of adding two roots of a fifth degree polynomial. Based on Martinelli's polynomial, it was possible to demonstrate Milanez's formula.
After finding the fifth degree resolvent equation of the sixth degree equation according to Milanez's relation, it is possible to find the roots of the sixth degree polynomial.
It is concluded that a formula that can solve specific cases of the sixth degree equation was successfully achieved.


\begin{thebibliography}{9}

\bibitem{zoladek} 

\begin{flushleft} ZOLADEK, H. The topological proof of Abel-Ruffini theorem. \textit{Topological Methods in Nonlinear Analysis}.v. 16, n.2, pp. 253-265, 2000. \\ \end{flushleft}
\bibitem{nickalls} 
\begin{flushleft}NICKALLS, R. W. D. A new approach to solving the cubic: \textit{Cardan's solution revealed}. The mathematical Gazette, v. 77, n. 480, pp. 354-359, 1993.\\ \end{flushleft}
\bibitem{martinelli} 
\begin{flushleft}ANDRADE, R. J. M. B. Resolução de uma equação do quinto grau. C.Q.D - Revista Eletrônica Paulista de Matemática, v. 16, pp. 196-205, dez. 2019. 
Retrieved from \url{https://www.fc.unesp.br/#!/departamentos/matematica/revista-cqd/}. \end{flushleft}
\end{thebibliography}
\end{document}